\documentclass[10pt]{article}
\usepackage{amssymb}
\usepackage{amsmath}
\usepackage[dvips]{graphics}
\usepackage{epsfig}
\begin{document}
\title{Interpolating the Derivatives of the Gamma Function}

\author{\textbf{Vassilis G. Papanicolaou}
\\\\
Department of Mathematics
\\
National Technical University of Athens,
\\
Zografou Campus, 157 80, Athens, GREECE
\\
{\tt papanico@math.ntua.gr}}
\maketitle

\begin{abstract}
We consider a function $G(\lambda, z)$, entire in $\lambda$, which interpolates the derivatives of the Gamma function in the sense that
$G(m, z) = \Gamma^{(m)}(z)$ for any integer $m \geq 0$ and we calculate the asymptotics of $G(\lambda, z)$ as $\lambda \to +\infty$.

\end{abstract}

\\\\
\textbf{2010 AMS Mathematics Classification.} 33B15; 41A60; 33B20.

\section{Introduction---the function $G(\lambda, z)$}
Starting with the familiar integral representation of the Gamma function
\begin{equation}
\Gamma(z) = \int_0^{\infty} e^{-t} t^{z-1} dt,
\qquad
\Re(z) > 0,
\label{A1}
\end{equation}
and then taking derivatives, one obtains
\begin{equation}
\Gamma^{(m)}(z) = \int_0^{\infty} e^{-t} (\ln t)^m t^{z-1} dt,
\qquad
m = 0, 1, 2, \ldots,
\quad
\Re(z) > 0.
\label{A2}
\end{equation}
The goal of this article is to construct and study an analytic function $G$ of two variables $\lambda$ and $z$ which, for $\lambda = m$ becomes
the $m$-th derivative $\Gamma^{(m)}(z)$ of the Gamma function.

To construct such a function we simply set
\begin{equation}
G(\lambda, z) := \int_{C_r} e^{-\zeta} (\ln \zeta)^{\lambda} \, \zeta^{z-1} d\zeta,
\qquad
\quad
\Re(z) > 0,
\label{A3b}
\end{equation}
where:

(i) The contour $C_r$ consists of the intervals $[0, 1-r]$ and $[1+r, \infty)$ of the positive real semiaxis ($0 < r < 1$) connected
with the semicircle $\{\zeta \in \mathbb{C} \,:\, \zeta = 1 + r e^{i\theta}, \  0 < \theta < \pi\}$, lying in the upper half-plane. The orientation of $C_r$ is $0 \to \infty$.

(ii) $w = \ln\zeta$ denotes the branch of the logarithm which maps the open upper half-plane $\Im(\zeta) > 0$ conformally onto the open strip
$0 < \Im(w) < \pi$. In particular, $\ln1 = 0$ and $\ln(-1) = \pi i$.

(iii) $(\ln \zeta)^{\lambda} = e^{\lambda \ln\ln \zeta}$ and $\zeta^{z-1} = e^{(z-1)\ln \zeta}$. In particular $1^{\lambda} = 1^{z-1} = 1$.

It is clear that the value of the integral in \eqref{A3b} is independent of $r \in (0, 1)$ and hence it defines a unique function $G(\lambda, z)$,
which is entire in $\lambda$. Moreover, it is not hard to see that, as a function of $\lambda$, $G(\lambda, z)$ is of order $1$ and of maximal (i.e.
infinite) type for any $z$ with $\Re(z) > 0$.

If $\Re(\lambda) > -1$ (and $\Re(z) > 0$), we can let $r \to 0^+$ in \eqref{A3b} and arrive at the formula
\begin{equation}
G(\lambda, z) = e^{\pi i \lambda} \int_0^1 e^{-t} (-\ln t)^{\lambda} \, t^{z-1} dt
+ \int_1^{\infty} e^{-t} (\ln t)^{\lambda} \, t^{z-1} dt,
\label{A3}
\end{equation}
since for the negative values of $\ln t$ we have
\begin{equation}
(\ln t)^{\lambda} = e^{\lambda \ln\ln t} = e^{\lambda [\ln(-\ln t) + \pi i]} = e^{\pi i \lambda} (-\ln t)^{\lambda}.
\label{A4}
\end{equation}
Formula \eqref{A3} can be also written as
\begin{equation}
G(\lambda, z) = \int_0^{\infty} e^{-t} (\ln t)^{\lambda} \, t^{z-1} dt,
\qquad
\Re(\lambda) > -1,
\quad
\Re(z) > 0,
\label{A3a}
\end{equation}
from which we can see immediately that $G(m, z) = \Gamma^{(m)}(z)$ for $m = 0, 1, \ldots$ and $\Re(z) > 0$.

\section{Some basic properties of $G(\lambda, z)$}
Before we start our analysis let us introduce the notation
\begin{equation}
\mathbb{C}_{\,\xi} := \{z \in \mathbb{C} \,:\, \Re(z) > \xi\},
\label{A4a}
\end{equation}
where $\xi$ is a real number.

\medskip

As we have, essentially, already observed, it follows from \eqref{A3b} that $G(\lambda, z)$ is analytic in
$\mathbb{C} \times \,\mathbb{C}_{\,0}$.

We continue with some basic properties of $G(\lambda, z)$.

\medskip

\textbf{Property 1.} For any $\lambda \in \mathbb{C}_{-1}$ we have
\begin{equation}
L(\lambda) :=
%\lim_{z \to 0 \atop z \in \mathbb{C}_{\,0}}
\lim_{\substack{z \to \, 0\\z \in \mathbb{C}_{\,0}}}
z^{\lambda+1} G(\lambda, z) = e^{\pi i \lambda} \Gamma(\lambda+1),
\label{A5}
\end{equation}
where $z^{\lambda+1} = e^{(\lambda+1) \ln z}$, $\ln z$ being continuous in $\mathbb{C}_{\,0}$ (in particular $1^{\lambda+1} = 1$).

\smallskip

\textit{Proof}. In view of \eqref{A3a} we have 
\begin{equation}
G(\lambda, z) = \int_0^{\varepsilon} e^{-t} (\ln t)^{\lambda} \, t^{z-1} dt +\int_{\varepsilon}^{\infty} e^{-t} (\ln t)^{\lambda} \, t^{z-1} dt,
\quad
\Re(\lambda) > -1,
\ 
\Re(z) > 0,
\label{A3aa}
\end{equation}
for any fixed $\varepsilon > 0$. Now, the integral
\begin{equation}
\int_{\varepsilon}^{\infty} e^{-t} (\ln t)^{\lambda} \, t^{z-1} dt
\label{A6}
\end{equation}
is entire in $z$ and, in particular, finite at $z=0$; thus it does not contribute to the value of the limit $L(\lambda)$. Consequently,
for the limit in \eqref{A5} we have
\begin{equation}
L(\lambda) = \lim_{\substack{z \to \, 0\\z \in \mathbb{C}_{\,0}}} z^{\lambda+1} \int_0^{\varepsilon} e^{-t} (\ln t)^{\lambda} \, t^{z-1} dt,
\label{A7}
\end{equation}
independently of $\varepsilon$, as long as $\varepsilon > 0$. Since for $t \in [0, \varepsilon]$ the quantity $e^{-t}$ can be made
as close to $1$ as we wish by choosing $\varepsilon$ sufficiently small, it follows that
\begin{equation}
L(\lambda) = \lim_{\substack{z \to \, 0\\z \in \mathbb{C}_{\,0}}} z^{\lambda+1} \int_0^{\varepsilon} (\ln t)^{\lambda} \, t^{z-1} dt.
\label{A8}
\end{equation}
But the limit in \eqref{A8} is, also, independent of $\varepsilon$, as long as $\varepsilon > 0$. Therefore, we can take $\varepsilon = 1$ in \eqref{A8} and get
\begin{align}
L(\lambda) &= \lim_{\substack{z \to \, 0\\z \in \mathbb{C}_{\,0}}} z^{\lambda+1} \int_0^1 (\ln t)^{\lambda} \, t^{z-1} dt
\nonumber
\\
&= e^{\pi i \lambda} \lim_{\substack{z \to \, 0\\z \in \mathbb{C}_{\,0}}} z^{\lambda+1} \int_0^1 (-\ln t)^{\lambda} \, t^{z-1} dt
\nonumber
\\
&= e^{\pi i \lambda} \lim_{\substack{z \to \, 0\\z \in \mathbb{C}_{\,0}}} z^{\lambda+1} \int_0^{\infty} e^{-z \tau} \tau^{\lambda} d\tau
\label{A9a}
\end{align}
(the last integral is obtained by substituting $\tau = -\ln t$ in the previous integral).

Finally, since for real and positive $z$, the substitution $t = z \tau$ yields
\begin{equation}
z^{\lambda+1} \int_0^{\infty} e^{-z \tau} \tau^{\lambda} d\tau = \Gamma(\lambda+1),
\qquad
\lambda \in \mathbb{C}_{-1},
\label{A9}
\end{equation}
it follows by analytic continuation that the equality \eqref{A9} holds for any $z \in \mathbb{C}_{\,0}$.

\hfill $\blacksquare$

\medskip

\textbf{Property 2.} If $\partial_z$ denotes the derivative operator with respect to $z$, then
\begin{equation}
\partial_z G(\lambda, z) = G(\lambda+1, z),
\qquad
\lambda \in \mathbb{C},
\quad
z \in \mathbb{C}_{\,0}.
\label{A10}
\end{equation}

\medskip

Property 2 follows by differentiating \eqref{A3b} with respect to $z$.

\medskip

\textbf{Property 3.} The function $G(\lambda, z)$ satisfies the functional equation
\begin{equation}
G(\lambda, z+1) = z G(\lambda, z) + \lambda G(\lambda-1, z),
\qquad
\lambda \in \mathbb{C},
\quad
z \in \mathbb{C}_{\,0}.
\label{A11}
\end{equation}

\smallskip

\textit{Proof}. For $\Re(\lambda) > 1$ and $\Re(z) > 0$ integration by parts yields
\begin{align}
z \int_0^{\infty} e^{-t} (\ln t)^{\lambda} \, t^{z-1} dt &= e^{-t} (\ln t)^{\lambda} \, t^z \Big|_{t=0}^{\infty}
+ \int_0^{\infty} \left[e^{-t} (\ln t)^{\lambda} \, t^z - \lambda e^{-t} (\ln t)^{\lambda-1} t^{z-1}\right] dt
\nonumber
\\
&= \int_0^{\infty} e^{-t} (\ln t)^{\lambda} \, t^z dt - \lambda \int_0^{\infty} e^{-t} (\ln t)^{\lambda-1} t^{z-1} dt
\nonumber
\\
&=G(\lambda, z+1) - \lambda G(\lambda-1, z).
%\label{A11a}
\nonumber
\end{align}
Thus, $G(\lambda, z)$ satisfies \eqref{A11} for $\Re(\lambda) > 1$ and, consequently, for all $\lambda \in \mathbb{C}$ by analytic continuation.
\hfill $\blacksquare$

\medskip

We can now combine Property 2 with Property 3 and get
\begin{equation}
G(\lambda, z+1) = z \partial_z G(\lambda-1, z) + \lambda G(\lambda-1, z),
\qquad
\lambda \in \mathbb{C}, \quad z \in \mathbb{C}_{\,0}.
\label{A12}
\end{equation}
Multiplying \eqref{A12} by $z^{\lambda-1}$ yields
\begin{equation}
z^{\lambda-1} G(\lambda, z+1) =  \partial_z \left[z^{\lambda} G(\lambda-1, z)\right],
\qquad
\lambda \in \mathbb{C}, \quad z \in \mathbb{C}_{\,0},
\label{A13}
\end{equation}
so that, by integrating \eqref{A13} from $z_0$ to $z$ we obtain
\begin{equation}
z^{\lambda} G(\lambda-1, z) - z_0^{\lambda} G(\lambda-1, z_0) = \int_{z_0}^z \zeta^{\lambda-1} G(\lambda, \zeta+1) \, d\zeta,
\qquad
\lambda \in \mathbb{C},
\label{A14}
\end{equation}
where the integral is taken over a rectifiable arc joining $z_0$ and $z$, and lying entirely in the open right half-plane $\mathbb{C}_{\,0}$.

Incidentally, if we restrict $\lambda \in \mathbb{C}_{\,0}$, then we can let $z_0 \to 0$ in \eqref{A14} and invoke Property 1 in order to get
\begin{equation}
z^{\lambda} G(\lambda-1, z) + e^{\pi i \lambda} \Gamma(\lambda) = \int_0^z \zeta^{\lambda-1} G(\lambda, \zeta+1) \, d\zeta,
\qquad
\lambda, z \in \mathbb{C}_{\,0}.
\label{A15}
\end{equation}
Noticing that, for any $\lambda \in \mathbb{C}$, the integral in the right-hand side of \eqref{A14} defines a (multivalued) analytic function of $z$
in the half-plane $\mathbb{C}_{-1}$ we can conclude that $z^{\lambda} G(\lambda-1, z)$ and, consequently,
$G(\lambda-1, z)$ have a meromorphic extension for $z \in \mathbb{C}_{-1}$, with a branch point at $z=0$. It is convenient (and harmless) to
denote by $G(\lambda-1, z)$ too the meromorphic extension of $G(\lambda-1, z)$.

It, then, follows that for a fixed $z_0 \in \mathbb{C}_{\,0}$ formula \eqref{A14} remains valid for
$z \in \mathbb{C}_{-1}$ (as long as the integral is taken over an arc joining $z_0$ and $z$, and lying entirely in $\mathbb{C}_{-1}$).
Therefore, by the same argument $G(\lambda-1, z)$ has a (multivalued) meromorphic extension for $z \in \mathbb{C}_{-2}$. As the dummy variable
$\zeta$ of the integral in \eqref{A14} takes the value $-1$, the function $G(\lambda, \zeta+1)$ appearing in the integrand is evaluated at its
branch point $0$. Hence $z=-1$ is a branch point of the integral and, hence, of $G(\lambda-1, z)$.
Keeping arguing in the same manner, we can conclude that $G(\lambda-1, z)$ has a (multivalued) meromorphic extension in $\mathbb{C}_{-n}$ for any
$n = 1, 2, \ldots$, with branch points at $z = 0, 1, \ldots, n-1$. And since $n$ is arbitrary, we can finally  conclude that $G(\lambda-1, z)$ has a
(multivalued) meromorphic extension in $\mathbb{C}$ with branch points at the nonpositive integers. For the integral values of $\lambda$,
namely for $\lambda = m$, $m = 1, 2, \ldots$, we know that $G(\lambda-1, z) = \Gamma^{(m-1)}(z)$ has poles of order $m$ at $z = -n$, where
$n = 0, 1, 2, \ldots\,$.

From the above discussion it also follows that the limit in \eqref{A5} of Proposition 1 remains valid as long as $z$ approaches zero in a continuous
fashion, without ever crossing the negative real semiaxis.

Let us summarize the above observations in the following property.

\medskip

\textbf{Property 4.} For any complex $\lambda$, the function $G(\lambda, z)$, viewed as a function of $z$, has a (multivalued) meromorphic extension in $\mathbb{C}$ with branch points at $z = 0, -1, -2, \ldots\,$. Furthermore, Property 1 can be strengthen as
\begin{equation}
\lim_{\substack{z \to \, 0\\z \in \gamma}}
z^{\lambda+1} G(\lambda, z) = e^{\pi i \lambda} \Gamma(\lambda+1),
\label{A16}
\end{equation}
where $\gamma$ is a continuous curve which does not cross the negative real semiaxis, while $z^{\lambda+1}$ is defined so that it is continuous in $z$ and $1^{\lambda+1} = 1$.

\medskip

Finally, from \eqref{A14} it follows easily the following property:

\medskip

\textbf{Property 5.} For any $z \ne 0, -1, -2, \ldots$, the function $G(\lambda, z)$ is entire in $\lambda$.

\section{Asymptotics of $G(\lambda, z)$ as $\lambda \to +\infty$}
For typographical convenience we write \eqref{A3} as
\begin{equation}
G(\lambda, z) = G_0(\lambda, z) + G_1(\lambda, z),
\label{B1}
\end{equation}
where
\begin{equation}
G_0(\lambda, z) := e^{\pi i \lambda} \int_0^1 e^{-t} (-\ln t)^{\lambda} \, t^{z-1} dt,
\qquad
\Re(\lambda) > -1,
\quad
\Re(z) > 0,
\label{B2a}
\end{equation}
and
\begin{equation}
G_1(\lambda, z) := \int_1^{\infty} e^{-(t - \lambda \ln\ln t)}  t^{z-1} dt,
\qquad
\Re(\lambda) > -1,
\quad
\Re(z) > 0
\label{B2b}
\end{equation}
(in fact, $G_1(\lambda, z)$ is entire in $z$).

\subsection{The expansion of $G_0(\lambda, z)$}
By expanding $e^{-t}$ in \eqref{B2a} as a power series we obtain
\begin{equation}
G_0(\lambda, z) = e^{\pi i \lambda} \sum_{n=0}^{\infty} \frac{(-1)^n}{n!} \int_0^1 (-\ln t)^{\lambda} \, t^{n+z-1} dt,
\qquad
\Re(\lambda) > -1,
\quad
\Re(z) > 0.
\label{B3}
\end{equation}
Now, as in \eqref{A9a}-\eqref{A9},
\begin{equation}
\int_0^1 (-\ln t)^{\lambda} \, t^{n+z-1} dt = \int_0^{\infty} e^{-(n+z) \tau} \tau^{\lambda} d\tau = \frac{\Gamma(\lambda+1)}{(n+z)^{\lambda+1}}.
\label{B4}
\end{equation}
Thus, \eqref{B3} becomes
\begin{equation}
G_0(\lambda, z) = e^{\pi i \lambda} \Gamma(\lambda+1) \sum_{n=0}^{\infty} \frac{(-1)^n}{n! \, (n+z)^{\lambda+1}},
\qquad
\Re(\lambda) > -1,
\quad
\Re(z) > 0,
\label{B5}
\end{equation}
and the series converges very rapidly. Since for $\Re(z) > 0$
\begin{equation*}
\frac{1}{(n+1+z)^{\lambda+1}} \ll \frac{1}{(n +z)^{\lambda+1}}
\qquad \text{as }\;
\lambda \to +\infty,
%\label{B5a}
\end{equation*}
for any $n = 0, 1, 2, \ldots$,
the expansion in \eqref{B5} describes completely the behavior of $G_0(\lambda, z)$  as $\lambda \to +\infty$. In particular \eqref{B5} implies that
for any fixed $z \in \mathbb{C}_0$ we have
\begin{equation}
G_0(\lambda, z) = e^{\pi i \lambda} \frac{\Gamma(\lambda+1)}{z^{\lambda+1}} \left[1 + O\left(\left(\frac{z}{z+1}\right)^{\lambda+1}\right)\right],
\qquad
\lambda \to +\infty.
\label{B6}
\end{equation}

\subsection{Asymptotics of $G_1(\lambda, z)$ as $\lambda \to +\infty$}
Let us, first, introduce some notation. We set
\begin{equation}
\psi(t) := t\ln t,
\qquad
t \in [1, \infty)
\label{B7a}
\end{equation}
and
\begin{equation}
\omega(\lambda) := \psi^{-1}(\lambda),
\qquad
\lambda \in [0, \infty).
\label{B7b}
\end{equation}
The function $\omega(\lambda)$ is well defined, smooth, and strictly increasing on $[0, \infty)$, with $\omega(0) = 1$, since $\psi(t)$ is smooth and strictly increasing on $[1, \infty)$, with $\psi(1) = 0$. Furthermore, it is easy to see that
\begin{equation}
\omega(\lambda) = \frac{\lambda}{\ln\lambda}\left[1 + O\left(\frac{\ln\ln\lambda}{\ln\lambda}\right)\right],
\qquad
\lambda \to +\infty.
\label{B7c}
\end{equation}

\medskip

\textbf{Lemma 1.} For any fixed $z \in \mathbb{C}_0$ and for any fixed $\alpha$ such that
\begin{equation}
0 < \alpha < \frac{1}{2},
\label{B16}
\end{equation}
the function $G_1(\lambda, z)$ of \eqref{B2b} satisfies
\begin{equation}
G_1(\lambda, z)
= \frac{\sqrt{\pi} \, e^{\lambda \ln\ln \omega(\lambda)} e^{-\omega(\lambda)}}{\sqrt{A}} \, \omega(\lambda)^z
\left[1 + O\left(\frac{1}{\lambda^{\alpha}}\right)\right],
\qquad
\lambda \to +\infty,
\label{B8}
\end{equation}
where $\omega(\lambda)$ is given by \eqref{B7a}-\eqref{B7b} and
\begin{equation}
A := \lambda \frac{1 + \ln \omega(\lambda)}{2\ln^2 \omega(\lambda)}
\label{SA4b}
\end{equation}
(thus, from \eqref{B7c} we have that $A \sim \lambda/(2 \ln \lambda)$ as $\lambda \to +\infty$).

\smallskip

\textit{Proof}. The integral of \eqref{B2b} can be viewed as a ``Laplace integral" with large parameter $\lambda$ \cite{B-O}. Thus,
we will try to apply the so-called Laplace method for integrals.

We begin by looking for the value of $t$ which minimizes the function
\begin{equation}
h(t) := t - \lambda \ln\ln t,
\qquad
t \in [1, \infty),
\label{B10}
\end{equation}
appearing in the exponent of the integrand in \eqref{B2b}. Since
\begin{equation}
h'(t) = 1 - \frac{\lambda}{t\ln t}
\label{B10a}
\end{equation}
we have that $h(t)$ has a unique minimum attained when
\begin{equation*}
t\ln t = \lambda
%\label{B10b}
\end{equation*}
namely when (recall \eqref{B7a}-\eqref{B7b})
\begin{equation}
t = \omega(\lambda)
\label{B11}
\end{equation}
and, since $\omega(\lambda) > 1$ for $\lambda > 0$, this minimum is attained in the interior of the interval $[1, \infty)$, for every $\lambda > 0$.

From the theory of Laplace integrals we know that the main contribution to the value of the integral in \eqref{B2b}, as $\lambda$ gets large,
comes from the values of $t$ around $\omega(\lambda)$. In order to avoid the dependence in $\lambda$ we make the substitution $t = \omega(\lambda)x$.
Then, formula \eqref{B2b} becomes
\begin{equation}
G_1(\lambda, z) = \omega(\lambda)^z
\int_{1/\omega(\lambda)}^{\infty} e^{-p(x)} x^{z-1}dx,
\label{B12}
\end{equation}
where for typographical convenience we have set
\begin{equation}
p(x) := h\big(\omega(\lambda)x\big) = \omega(\lambda)x - \lambda \ln(\ln \omega(\lambda) + \ln x),
\qquad
x \in \left[1/\omega(\lambda), \, \infty\right).
\label{B13}
\end{equation}
The minimum of $p(x)$ on $[1/\omega(\lambda), \infty)$ is attained at $x=1$. Since $p'(1) = 0$
the Taylor expansion of $p(x)$ (with remainder) about $x=1$ up to the cubic term is
\begin{equation}
p(x) = p(1) + \frac{p''(1)}{2} (x-1)^2 + \frac{p'''(c)}{6} (x-1)^3,
\label{B14}
\end{equation}
where $c$ is some number between $1$ and $x$.

From \eqref{B13} we have
\begin{equation}
p(1) = \omega(\lambda) - \lambda \ln\ln \omega(\lambda),
\qquad\qquad
p''(1) = \lambda \frac{1 + \ln \omega(\lambda)}{\ln^2 \omega(\lambda)},
\label{B14a}
\end{equation}
and
\begin{equation}
p'''(c) = -\lambda \left\{\frac{2[\ln \omega(\lambda) + \ln c]^2 + 3[\ln \omega(\lambda) + \ln c] +2}{c^3 [\ln \omega(\lambda) + \ln c]^3} \right\}.
\label{B15}
\end{equation}
In particular, formula \eqref{B15}, combined with \eqref{B7c}, imply that
\begin{equation}
p'''(c) = O\left(\frac{\lambda}{\ln\lambda}\right)
\qquad \text{as }\; \lambda \to +\infty,
\label{B15a}
\end{equation}
provided $\delta \leq x \leq 1/\delta$ for some fixed $\delta \in (0, 1)$. Moreover, by \eqref{B14} and \eqref{B14a} we can see that the
quadratic Taylor polynomial of $p(x)$ about $x=1$ is
\begin{equation}
P(x; \lambda) := \omega(\lambda) - \lambda \ln\ln \omega(\lambda) + \lambda \frac{1 + \ln \omega(\lambda)}{2\ln^2 \omega(\lambda)} (x-1)^2
\label{B18}
\end{equation}
and it follows easily from \eqref{B14}, \eqref{B14a}, \eqref{B15}, \eqref{B18}, and \eqref{B16} that for all $\lambda$ sufficiently large we must
have
\begin{equation}
\left|p(x) - P(x; \lambda)\right|
\leq \frac{1}{\lambda^{3\alpha - 1}},
\qquad
x \in I_{\lambda},
\label{B17a}
\end{equation}
where $I_{\lambda}$ is the interval
\begin{equation}
I_{\lambda} := \left[1 - \frac{1}{\lambda^{\alpha}}, 1 + \frac{1}{\lambda^{\alpha}}\right].
\label{B17b}
\end{equation}
From now on let us assume, without loss of generality, that $\alpha > 1/3$, so that $1/\lambda^{3\alpha - 1} \to 0$ as $\lambda \to +\infty$. Then,
from \eqref{B17a} it follows that
\begin{equation}
e^{-p(x)} = e^{-P(x; \lambda)} + R(x; \lambda),
\qquad
x \in I_{\lambda},
\label{B19a}
\end{equation}
where
\begin{equation}
|R(x; \lambda)| \leq \frac{2}{\lambda^{3\alpha - 1}} \, e^{-P(x; \lambda)},
\qquad
x \in I_{\lambda},
\label{B19b}
\end{equation}
for all sufficiently large $\lambda$.

Multiplying by $x^{z-1}$ and then taking integrals in \eqref{B19a} yields
\begin{equation}
\int_{1-\lambda^{-\alpha}}^{1+\lambda^{-\alpha}} e^{-p(x)} x^{z-1} dx
= \int_{1-\lambda^{-\alpha}}^{1+\lambda^{-\alpha}} e^{-P(x; \lambda)} x^{z-1} dx
+ \int_{1-\lambda^{-\alpha}}^{1+\lambda^{-\alpha}} R(x; \lambda) x^{z-1}dx.
\label{B20}
\end{equation}

Having formula \eqref{B20}, in order to complete the proof of the lemma, we need to verify the following three claims.

\smallskip

\textbf{Claim 1.} For the first integral in the right-hand side of \eqref{B20} we have
\begin{equation}
\int_{1-\lambda^{-\alpha}}^{1+\lambda^{-\alpha}} e^{-P(x; \lambda)} x^{z-1} dx
= \frac{\sqrt{\pi} \, e^{-\omega(\lambda) + \lambda \ln\ln \omega(\lambda)}}{\sqrt{A}} \left[1 + O\left(\frac{1}{\lambda^{\alpha}}\right)\right],
\qquad
\lambda \to +\infty.
\label{SA1}
\end{equation}

\smallskip

\textbf{Claim 2.} For the second integral in the right-hand side of \eqref{B20} we have
\begin{equation}
\int_{1-\lambda^{-\alpha}}^{1+\lambda^{-\alpha}} R(x; \lambda) x^{z-1}dx
= \frac{e^{-\omega(\lambda) + \lambda \ln\ln \omega(\lambda)}}{\sqrt{A}} \, O\left(\frac{1}{\lambda^{3\alpha-1}}\right),
\qquad
\lambda \to +\infty.
\label{SB1}
\end{equation}

\smallskip

\textbf{Claim 3.} For the integral in the left-hand side of \eqref{B20} we have
\begin{equation}
\int_{1-\lambda^{-\alpha}}^{1+\lambda^{-\alpha}} e^{-p(x)} x^{z-1}dx
= \int_{1/\omega(\lambda)}^{\infty} e^{-p(x)} x^{z-1} dx
+ \frac{e^{-\omega(\lambda) + \lambda \ln\ln \omega(\lambda)}}{\sqrt{A}} \, o\left(\frac{1}{\lambda^{\alpha}}\right)
\label{SC1}
\end{equation}
as $\lambda \to +\infty$.

\smallskip

\underline{\textit{Proof of Claim 1}}. In view of \eqref{B18} we have
\begin{equation}
\int_{1-\lambda^{-\alpha}}^{1+\lambda^{-\alpha}} e^{-P(x; \lambda)} x^{z-1} dx
= e^{-\omega(\lambda) + \lambda \ln\ln \omega(\lambda)}
\int_{-\lambda^{-\alpha}}^{\lambda^{-\alpha}} e^{-\lambda \frac{1 + \ln \omega(\lambda)}{2\ln^2 \omega(\lambda)} \, \xi^2} (1+\xi)^{z-1} d\xi.
\label{SA2}
\end{equation}
Now, for $|\xi| \leq \lambda^{-\alpha}$ we have (since $z$ is fixed)
\begin{equation}
(1+\xi)^{z-1} = 1 + O\left(\frac{1}{\lambda^{\alpha}}\right),
\qquad
\lambda \to +\infty.
\label{SA3}
\end{equation}
Thus, \eqref{SA2} implies
\begin{equation}
\int_{1-\lambda^{-\alpha}}^{1+\lambda^{-\alpha}} e^{-P(x; \lambda)} x^{z-1} dx
= e^{-\omega(\lambda) + \lambda \ln\ln \omega(\lambda)} \left[1 + O\left(\frac{1}{\lambda^{\alpha}}\right)\right]
\int_{-\lambda^{-\alpha}}^{\lambda^{-\alpha}} e^{-A \xi^2} d\xi,
\label{SA4a}
\end{equation}
where $A$ is given by \eqref{SA4b}.

Finally, since from \eqref{B16}, \eqref{B7c}, and \eqref{SA4b} we have that $\sqrt{A} \, \lambda^{-\alpha} \gg \lambda^{(1-2\alpha)/4} \to +\infty$,
while
\begin{equation}
\int_0^{\sqrt{A} \, \lambda^{-\alpha}} e^{-u^2} du
= \frac{\sqrt{\pi}}{2} + O\left(\frac{e^{-A \, \lambda^{-2\alpha}}}{\sqrt{A} \, \lambda^{-\alpha}}\right)
\label{SA6}
\end{equation}
as $\lambda \to +\infty$, formula \eqref{SA4a} implies \eqref{SA1}.

\smallskip

\underline{\textit{Proof of Claim 2}}. Using \eqref{B19b} we get
\begin{equation}
\left|\int_{1-\lambda^{-\alpha}}^{1+\lambda^{-\alpha}} R(x; \lambda) x^{z-1}dx\right|
\leq \frac{2}{\lambda^{3\alpha - 1}}
\int_{1-\lambda^{-\alpha}}^{1+\lambda^{-\alpha}} e^{-P(x; \lambda)} x^{\Re(z)-1}dx,
\label{B21}
\end{equation}
thus \eqref{SB1} follows by using \eqref{SA1} in \eqref{B21}.

\smallskip

\underline{\textit{Proof of Claim 3}}. We assume, without loss of generality, that $\lambda > 1$. Then
\begin{equation}
\left|\int_{1+\lambda^{-\alpha}}^{\infty} e^{-p(x)} x^{z-1}dx\right|
\leq \int_{1+\lambda^{-\alpha}}^2 e^{-p(x)} x^{\Re(z)-1}dx + \int_2^{\infty} e^{-p(x)} x^{\Re(z)-1}dx
\label{SC2}
\end{equation}
From \eqref{B13} we have that $p''(x) > 0$ for all $x \in \left[1/\omega(\lambda), \, \infty\right)$, while $p'(1) = 0$ by \eqref{B7a}-\eqref{B7b}. Thus, $p(x)$ is convex on $[1/\omega(\lambda), \, \infty)$ and increasing on $[1, \infty)$. Consequently, for the first integral in the right-hand
side of \eqref{SC2} we have
\begin{align}
\int_{1+\lambda^{-\alpha}}^2 e^{-p(x)} x^{\Re(z)-1}dx
&\leq 2^{\Re(z)} \int_{1+\lambda^{-\alpha}}^2 e^{-p(x)} dx
\nonumber
\\
&\leq 2^{\Re(z)} e^{-p(1+\lambda^{-\alpha})} \int_{1+\lambda^{-\alpha}}^2 e^{-p'(1+\lambda^{-\alpha}) (x - 1 - \lambda^{-\alpha})} dx
\label{SC2a}
\end{align}
(since $p(x) \geq p(1+\lambda^{-\alpha}) + p'(1+\lambda^{-\alpha}) (x - 1 - \lambda^{-\alpha})$ on $\left[1/\omega(\lambda), \, \infty\right)$ by
convexity).

Formula \eqref{SC2a} implies
\begin{equation}
\int_{1+\lambda^{-\alpha}}^2 e^{-p(x)} x^{\Re(z)-1}dx
\leq 2^{\Re(z)} \frac{e^{-p(1+\lambda^{-\alpha})}}{p'(1+\lambda^{-\alpha})}.
\label{SC4}
\end{equation}
By \eqref{B17a} we have
\begin{equation}
p(1+\lambda^{-\alpha})
= \omega(\lambda) - \lambda \ln\ln \omega(\lambda) + \frac{1 + \ln \omega(\lambda)}{2\ln^2 \omega(\lambda)} \lambda^{1-2\alpha}
+ O\left(\frac{1}{\lambda^{3\alpha - 1}}\right),
\label{SC5}
\end{equation}
thus
\begin{equation}
e^{-p(1+\lambda^{-\alpha})}
= e^{-\omega(\lambda) + \lambda \ln\ln \omega(\lambda)} \, e^{-\frac{1 + \ln \omega(\lambda)}{2\ln^2 \omega(\lambda)} \lambda^{1-2\alpha}}
\left[1 + O\left(\frac{1}{\lambda^{3\alpha - 1}}\right)\right].
\label{SC6}
\end{equation}
Furthermore, from the Taylor expansion with remainder of $p'(x)$ about $x=1$, namely
\begin{equation*}
p'(x) = p'(1) + p''(1) (x-1) + \frac{p'''(c)}{2} (x-1)^2,
\qquad
x \in I_{\lambda},
%\label{B14aa}
\end{equation*}
together with \eqref{B14a}, \eqref{B15a}, and the fact that $p'(1) = 0$, we obtain
\begin{equation}
p'(1+\lambda^{-\alpha}) = \frac{\lambda^{1-\alpha}}{\ln \omega(\lambda)} \left[1 + O\left(\lambda^{1-2\alpha}\right)\right].
\label{SC7}
\end{equation}
Hence, by using \eqref{SC6} and \eqref{SC7} in \eqref{SC4} (and recalling \eqref{B7c} and the fact that $1-2\alpha > 0$) it follows that
\begin{equation}
\int_{1+\lambda^{-\alpha}}^2 e^{-p(x)} x^{\Re(z)-1}dx
= \frac{e^{-\omega(\lambda) + \lambda \ln\ln \omega(\lambda)}}{\sqrt{A}} \, o\left(\frac{1}{\lambda^{\alpha}}\right),
\qquad
\lambda \to +\infty.
\label{SC8}
\end{equation}
As for the second integral in the right-hand side of \eqref{SC2}, it is easy to see that it is much smaller than the
bound given by \eqref{SC8} for the first integral in the right-hand side of \eqref{SC2}. We can see that, e.g., by writing
\begin{equation}
\int_2^{\infty} e^{-p(x)} x^{\Re(z)-1}dx = \int_2^{\infty} e^{-[p(x) - p(x - 3/2)]} e^{-p(x - 3/2)} x^{\Re(z)-1}dx
\label{SC9}
\end{equation}
and observing that the function $e^{-p(x - 3/2)} x^{\Re(z)-1}$ is bounded on $[2, \infty)$ uniformly in $\lambda$ (for, say, $\lambda > 2$).

Finally, in the same manner we can show that
\begin{equation}
\int_{1/\omega(\lambda)}^{1-\lambda^{-\alpha}} e^{-p(x)} x^{z-1}dx
= \frac{e^{-\omega(\lambda) + \lambda \ln\ln \omega(\lambda)}}{\sqrt{A}} \, o\left(\frac{1}{\lambda^{\alpha}}\right),
\qquad
\lambda \to +\infty.
\label{SC10}
\end{equation}
Therefore, \eqref{SC1} follows immediately by using \eqref{SC8} and \eqref{SC10} in \eqref{SC2}.
\hfill $\blacksquare$

\medskip

If we use the estimate \eqref{B7c} in \eqref{B8} we get the following corollary

\medskip

\textbf{Corollary 1.} For the function $G_1(\lambda, z)$ of \eqref{B2b} we have
\begin{equation}
G_1(\lambda, z)
\sim \sqrt{\frac{\pi}{2}} \,\, \big[\ln \omega(\lambda)\big]^{\lambda} e^{-\omega(\lambda)} \left(\frac{\lambda}{\ln\lambda}\right)^{z-(1/2)}
\qquad \text{as }\;
\lambda \to +\infty,
\label{C1}
\end{equation}
where $\omega(\lambda)$ is given by \eqref{B7a}-\eqref{B7b}.

\medskip

\textbf{Remark 1.} For the function $G_1(\lambda, z)$ of \eqref{B2b} we have
\begin{equation}
G_1(m, z) = \Gamma^{(m)}(z, 1),
\qquad
m = 0, 1, \ldots,
\label{C2}
\end{equation}
where
\begin{equation}
\Gamma(z, 1) := \int_1^{\infty} e^{-t} t^{z-1} dt,
\qquad
z \in \mathbb{C},
\label{C3}
\end{equation}
is the so-called upper incomplete Gamma function. Since $\Gamma(z, 1)$ is entire in $z$, the radius of convergence of its Taylor series about any
$z_0$ is infinite, and this agrees with the asymptotics of $G_1(\lambda, z)$ for large $\lambda$, as given by \eqref{C1}.

As for the function $G_0(\lambda, z)$ of \eqref{B2a} we have
\begin{equation}
G_0(m, z) = \gamma^{(m)}(z, 1),
\qquad
m = 0, 1, \ldots,
\label{C4}
\end{equation}
where
\begin{equation}
\gamma(z, 1) := \int_0^1 e^{-t} t^{z-1} dt,
\qquad
\Re(z) > 0,
\label{C5}
\end{equation}
is the so-called lower incomplete Gamma function. Since $\gamma(z, 1)$ has a simple pole at $z = 0$, the radius of convergence of its Taylor series
about any $z_0$ with $\Re(z) > 0$ is $|z_0|$, and this agrees with the asymptotics of $G_0(\lambda, z)$ for large $\lambda$, as given by \eqref{B6}.
Thus $G_0(\lambda, z)$ grows much faster than $G_1(\lambda, z)$ as $\lambda \to \infty$. This is, also, reflected in the corollary that follows.

\medskip

\textbf{Corollary 2.} For the function $G(\lambda, z)$ of \eqref{A3b} we have
\begin{align}
G(\lambda, z) =& \, e^{\pi i \lambda} \Gamma(\lambda+1) \sum_{n=0}^{\infty} \frac{(-1)^n}{n! \, (n+z)^{\lambda+1}}
\nonumber
\\
&+ \sqrt{\frac{\pi}{2}} \,\, \big[\ln \omega(\lambda)\big]^{\lambda} e^{-\omega(\lambda)} \left(\frac{\lambda}{\ln\lambda}\right)^{z-(1/2)}
[1 + o(1)],
\qquad
\lambda \to +\infty,
\label{C6}
\end{align}
where $z$ is fixed with $\Re(z) > 0$ and $\omega(\lambda)$ is given by \eqref{B7a}-\eqref{B7b}.

\medskip

Corollary 2 follows immediately from \eqref{B5} and Corollary 1. Since $G(m, z) = \Gamma^{(m)}(z)$ we obtain immediately from \eqref{C6} the behavior of $\Gamma^{(m)}(z)$ as $m \to \infty$.

\medskip

\textbf{Corollary 3.} For $\Re(z) > 0$ we have
\begin{align}
\Gamma^{(m)}(z) =& \, (-1)^m m! \sum_{n=0}^{\infty} \frac{(-1)^n}{n! \, (n+z)^{m+1}}
\nonumber
\\
&+ \sqrt{\frac{\pi}{2}} \,\, \big[\ln \omega(m)\big]^m e^{-\omega(m)} \left(\frac{m}{\ln m}\right)^{z-(1/2)}
[1 + o(1)],
\qquad
\lambda \to +\infty,
\label{C7}
\end{align}
where $\omega(\cdot)$ is given by \eqref{B7a}-\eqref{B7b}. In particular,
\begin{equation}
\Gamma^{(m)}(z) \sim \frac{(-1)^m m!}{z^{m+1}},
\qquad
\lambda \to +\infty.
\label{C7a}
\end{equation}

\subsection{Examples}

\textbf{1.} If we set $z=1$ in \eqref{C7}, we obtain
\begin{align}
\Gamma^{(m)}(1) =& \, (-1)^m m! \sum_{n=0}^{\infty} \frac{(-1)^n}{n! \, (n+1)^{m+1}}
\nonumber
\\
&+ \sqrt{\frac{\pi}{2}} \,\, \big[\ln \omega(m)\big]^m e^{-\omega(m)} \sqrt{\frac{m}{\ln m}} \, [1 + o(1)],
\qquad
\lambda \to +\infty,
\label{C8}
\end{align}
in particular, $\Gamma^{(m)}(1) \sim (-1)^m m!$.
Also, in view of \eqref{C2}, formula \eqref{C1} gives
\begin{equation}
\Gamma^{(m)}(1, 1)
\sim \sqrt{\frac{\pi}{2}} \,\, \big[\ln \omega(m)\big]^m e^{-\omega(m)} \sqrt{\frac{m}{\ln m}}
\qquad \text{as }\;
m \to \infty,
\label{C9}
\end{equation}
where $\omega(\cdot)$ is given by \eqref{B7a}-\eqref{B7b}.

\medskip

\textbf{2.} Suppose we want the asymptotics of $G(\lambda, z)$ as $\lambda \to \infty$, in the case where $z$ is a given complex number with
$\Re(z) < 0$. Then, we can employ formula \eqref{A14} or \eqref{A15}. For instance, for $z = -1/2$ formula \eqref{A15} gives
\begin{equation*}
G\left(\lambda-1, -\frac{1}{2}\right)
= -2^{\lambda} \Gamma(\lambda) - 2^{\lambda} e^{-\pi i \lambda} \int_{-1/2}^0 x^{\lambda-1} G(\lambda, x+1) \, dx,
%\label{C10a}
\end{equation*}
or
\begin{equation}
G\left(\lambda-1, -\frac{1}{2}\right)
= -2^{\lambda} \Gamma(\lambda) - e^{-\pi i \lambda} \int_0^1 (\xi-1)^{\lambda-1} G\left(\lambda, \frac{\xi+1}{2}\right) \, d\xi.
\label{C10}
\end{equation}
We can now use \eqref{C6} in \eqref{C10} and get (as $\lambda \to +\infty$)
\begin{align}
G\left(\lambda-1, -\frac{1}{2}\right) =& -2^{\lambda} \Gamma(\lambda) + e^{\pi i \lambda} 2^{\lambda+1} \Gamma(\lambda+1) J(\lambda)
\nonumber
\\
&+ \left[\sqrt{\frac{\pi}{2}} + o(1)\right] \big[\ln \omega(\lambda)\big]^{\lambda} e^{-\omega(\lambda)}
\int_0^1 (\xi-1)^{\lambda-1} \left(\frac{\lambda}{\ln\lambda}\right)^{\xi/2} d\xi 
\label{C11}
\end{align}
where we have set
\begin{equation}
J(\lambda) := \sum_{n=0}^{\infty} \frac{(-1)^n}{n!}
\int_0^1 \frac{(1-\xi)^{\lambda-1}}{(2n+1+\xi)^{\lambda+1}} \, d\xi.
\label{C12}
\end{equation}
By invoking the formula
\begin{equation}
-\frac{1}{(a+b)\lambda} \, \frac{d}{d\xi} \left[\frac{(a-\xi)^{\lambda}}{(b+\xi)^{\lambda}}\right]
= \frac{(a-\xi)^{\lambda-1}}{(b+\xi)^{\lambda+1}},
\label{C13}
\end{equation}
\eqref{C12} becomes
\begin{equation}
J(\lambda) = \frac{1}{2\lambda} \sum_{n=0}^{\infty} \frac{(-1)^n}{(n+1)!} \frac{1}{(2n+1)^{\lambda}}.
\label{C14}
\end{equation}
Also, it is not hard to show that
\begin{equation}
\int_0^1 (\xi-1)^{\lambda-1} \left(\frac{\lambda}{\ln\lambda}\right)^{\xi/2} d\xi \sim \frac{1}{\lambda},
\qquad
\lambda \to +\infty.
\label{C15}
\end{equation}
Using \eqref{C14} and \eqref{C15} in \eqref{C11} yields
\begin{align}
G\left(\lambda-1, -\frac{1}{2}\right) =& \left(e^{\pi i \lambda} - 1\right) 2^{\lambda} \Gamma(\lambda) + e^{\pi i \lambda} 2^{\lambda} \Gamma(\lambda)
\sum_{n=1}^{\infty} \frac{(-1)^n}{(n+1)!} \frac{1}{(2n+1)^{\lambda}}
\nonumber
\\
&+ \frac{\big[\ln \omega(\lambda)\big]^{\lambda}}{\lambda} \, e^{-\omega(\lambda)} \left[\sqrt{\frac{\pi}{2}} + o(1)\right],
\qquad
\lambda \to +\infty.
\label{C16}
\end{align}
By analytic continuation we have that $G(m, z) = \Gamma^{(m)}(z)$ for all $z \ne 0, -1, \ldots\,$. Hence, \eqref{C16} implies
\begin{align}
\Gamma^{(m)}\left(-\frac{1}{2}\right)
=& \, 2^{m+1} m! \left\{\left[(-1)^{m+1} - 1\right] - (-1)^m 
\sum_{n=1}^{\infty} \frac{(-1)^n}{(n+1)!} \frac{1}{(2n+1)^{m+1}}\right\}
\nonumber
\\
&+ \frac{\big[\ln \omega(m+1)\big]^{m+1}}{m+1} \, e^{-\omega(m+1)} \left[\sqrt{\frac{\pi}{2}} + o(1)\right],
\qquad
m \to \infty,
\label{C17}
\end{align}
in particular,
\begin{equation*}
\Gamma^{(2k)}\left(-\frac{1}{2}\right) \sim -2^{2k+2} (2k)!
\qquad \text{and} \qquad
\Gamma^{(2k+1)}\left(-\frac{1}{2}\right) \sim -\frac{2^{2k+1} (2k+1)!}{3^{2k+2}}
%\label{C17a}
\end{equation*}
as $k \to \infty$. The above asymptotic formulas are in agreement with the fact that for $z \in (-1, 0)$ the ``dominant" part of $\Gamma(z)$ is
\begin{equation*}
D_{\Gamma}(z) := \frac{1}{z} - \frac{1}{z+1}
%\label{C18}
\end{equation*}
for which we have $D_{\Gamma}^{(2k)}(-1/2) = -2^{2k+2} (2k)!$, while $D_{\Gamma}^{(2k+1)}(-1/2) = 0$. Let us also mention that from the functional
equation
\begin{equation}
\Gamma^{(m)}(z+1) = z \Gamma^{(m)}(1) + m \Gamma^{(m-1)}(z),
\qquad
z \in \mathbb{C},
\quad
m =0, 1, \ldots,
\label{C19}
\end{equation}
and the limits
\begin{equation}
\Gamma^{(2k)}(0^-) = \Gamma^{(2k)}(-1^+) = -\infty,
\qquad
\Gamma^{(2k+1)}(0^-) = -\Gamma^{(2k)}(-1^+) = -\infty,
\label{C20}
\end{equation}
it follows that $\Gamma^{(2k)}(x) < 0$ for $x \in (-1, 0)$ and hence that $\Gamma^{(2k+1)}(x)$ is decreasing on $(-1, 0)$. Consequently, for each
$k = 0, 1, \ldots$ the function $\Gamma^{(2k+1)}(z)$ has a unique zero, say $\eta_k$, in $(-1, 0)$. A natural question here is whether $\eta_k$,
$k \geq 0$, is a monotone sequence and, furthermore, whether $\eta_k \to -1/2$ as $k \to \infty$.

\medskip

\textbf{3.} Let us, finally, consider the case where $z$ is real and positive. Then formula \eqref{A2} implies that $\Gamma^{(2k)}(z) > 0$ for all
$k \geq 0$. Consequently, all odd derivatives $\Gamma^{(2k+1)}(z)$, $k \geq 0$, are increasing in $(0, \infty)$, with
$\Gamma^{(2k+1)}(0^+) = -\infty$ and $\Gamma^{(2k+1)}(+\infty) = +\infty$. Hence, for each $k \geq 0$ there is a unique $\zeta_k \in (0, \infty)$ such that
\begin{equation}
\Gamma^{(2k+1)}(\zeta_k) = 0.
\label{D1}
\end{equation}
Suppose that the sequence $\zeta_k$, $k \geq 0$, has a bounded subsequence, namely there is an $M > 0$ such that $\zeta_k < M$ for infinitely many
values of $k$. Then, we should have $\Gamma^{(2k+1)}(M) > 0$ for infinitely many values of $k$. But this is impossible, since formula \eqref{C7}, or
just \eqref{C7a}, implies that $\Gamma^{(2k+1)}(M) < 0$ for all sufficiently large $k$. Therefore $\zeta_k \to \infty$. A natural question here is
whether $\zeta_k$ is increasing.

\end{document}